\documentclass[11pt,notitlepage,twoside,a4paper]{amsart}
 \usepackage{amsfonts}
\usepackage[latin1]{inputenc}
\usepackage[cyr]{aeguill}
\usepackage{xspace}
\usepackage[polutonikogreek,english]{babel}

\usepackage{amsmath,amssymb,enumerate}

\usepackage{epsfig,fancyhdr,color}

\usepackage{amssymb}
\usepackage{amsmath,amsthm}
\usepackage{latexsym}
\usepackage{amscd}
\usepackage{psfrag}
\usepackage{graphicx}
\usepackage[latin1]{inputenc}
\usepackage[all]{xy}
\usepackage[mathcal]{eucal}

\definecolor{NoteColor}{rgb}{1,0,0}

\renewcommand{\textsc}{\textcolor{red}}
\newcommand*{\tg}[1]{\textgreek{#1}}

\newtheorem*{theorem 1}{\rm\bf Proposition 1}
\newtheorem*{theorem 2}{\rm\bf Proposition 2}

\theoremstyle{definition}

\theoremstyle{remark}

\def\interieur#1{\mathord{\mathop{\kern 0pt #1}\limits^\circ}}

\title[Convexity in Greek antiquity]{Convexity in Greek antiquity}

\author{Athanase Papadopoulos}
\address{Athanase Papadopoulos,  Universit{\'e} de Strasbourg and CNRS,
7 rue Ren\'e Descartes,
 67084 Strasbourg Cedex, France}
\email{athanase.papadopoulos@math.unistra.fr}

 \makeindex
\date{\today}

\begin{document}

 \begin{abstract}  
  We consider several appearances of the notion of convexity in Greek antiquity, more specifically in mathematics and optics, in the writings of Aristotle, and in art.
  
  The final version of this article will appear in the book \emph{Geometry in History}, ed. S. G. Dani and A. Papadopoulos, Springer Verlag, 2019.

       

  \end{abstract}  
    \maketitle
  \medskip

        \tableofcontents
    \vskip .3in

    \section{Introduction}
    
The mathematical idea of convexity was known in ancient Greece. It is present in the works on geometry and on optics of  Euclid\index{Euclid (c. 325--c. 270 BCE)} (c. 325--270 BCE),  Archimedes  (c. 287--212 BCE), Apollonius of Perga (c. 262--190 BCE), Heron of Alexandria (c. 10--70 CE), Ptolemy (c. 100--160 CE), and other mathematicians. This concept evolved slowly until the modern period where progress  was made by Kepler, Descartes and Euler, and convexity became gradually a property at the basis of several geometric results. For instance, Euler, in his memoir  \emph{De linea brevissima in superficie quacunque duo quaelibet puncta jungente} (Concerning the shortest line on any surface by which any two points can be joined together) (1732), gave the differential equation satisfied by a geodesic joining two points on a differentiable convex surface. 

     Around the beginning of the twentieth century, convexity acquired the status of a mathematical field, with works of Minkowski, Carathéodory, Steinitz, Fenchel, Jessen, Alexandrov, Busemann, Pogorelov and others.    
     
      My purpose in this note is to indicate some instances where the notion of convexity appears in the writings of the Greek mathematicians and philosophers of antiquity. The account is not chronological, because I wanted to start with convexity in the purely mathematical works, before talking about this notion in philosophy, architecture, etc.
     
     The present account of the history of convexity is very different from the existing papers on the subject. For sources concerning the modern period, I refer the reader to the paper by Fenchel \cite{Fenchel}.
     
       The final version of this article will appear in the book \emph{Geometry in History}, ed. S. G. Dani and A. Papadopoulos, Springer Verlag, 2019.

\section{Geometry}

Euclid\index{Euclid (c. 325--270 BCE)} uses convexity\index{convexity} in the \emph{Elements}\index{Euclid (c. 325--270 BCE)!\emph{Elements}}, although he does not give any precise definition of this notion.
In Proposition 8 of Book III, the words \emph{concave} and \emph{convex} sides of a  circumference appear, and Euclid\index{Euclid (c. 325--270 BCE)} regards them as understood. Propositions  36 and 37 of the same book also involve convexity:  Euclid talks about lines falling on the convex circumference. The constructibility of certain convex regular polygons is extensively studied in Book IV of the same work.  Books XI, XII and XIII are dedicated to the construction and properties of convex regular polyhedra. The word ``convex" is not used there to describe a property of these polyhedra, but Euclid relies extensively on the existence of a circumscribed sphere, which (in addition to the other properties that these polyhedra satisfy) implies that the polyhedra are convex. Such a sphere is also used by him for addressing the construction question: to construct the edge length of a face of a regular polyhedron in terms of the radius of the circumscribed sphere. 

There is a rather long passage on the construction and the properties of the regular convex polyhedra derived from their plane faces in Plato's\index{Plato (428--348 BCE)} \emph{Timaeus},\index{Plato (428--348 BCE)!\emph{Timaeus}}\footnote{\emph{Timaeus}, \cite{Cornford-Timaeus}, 53C-55C.} written about half-a-century before Euclid\index{Euclid (c. 325--270 BCE)}'s \emph{Elements}\index{Euclid (c. 325--270 BCE)!\emph{Elements}} appeared,  and it is commonly admitted that Plato learned this theory from Theaetetus,\index{Theaetetus (c. 417--c. 369 BCE)} a mathematician who was like him a student of the Pythagorean geometer Theodorus of Cyrene.\index{Theodorus of Cyrene (c. 465--398 BCE)} Actually, Plato\index{Plato (428--348 BCE)}\index{Plato (428--348 BCE)!\emph{Timaeus}} in the \emph{Timaeus}, was mainly interested in the construction of four out of the five regular polyhedra that he assigned to the four elements of nature, namely, the tetrahedron, the octahedron, the icosahedron and the cube. In particular, he shows how the faces of these polyhedra decompose into known (constructible) triangles, he computes angles between faces, etc. \cite[p. 210ff]{Cornford-Timaeus}.
The regular convex polyhedra were part of the teaching of the early Pythagoreans (cf. \cite{Iamblichus}). All this was discussed at length by several authors and commentators, see in particular Heath's notes in his edition of the \emph{Elements}\index{Euclid (c. 325--270 BCE)!\emph{Elements}} \cite{Euclid-Heath} and Cornford's comments in his edition of the \emph{Timaeus}.\index{Plato (428--348 BCE)!\emph{Timaeus}}

Apollonius,\index{Apollonius of Perga (c. 262--c. 190 BCE)}  in the \emph{Conics},\index{Apollonius of Perga (c. 262--c. 190 BCE)!\emph{Conics}} uses the notion of convexity,\index{convexity}  in particular in Book IV where he studies intersections of conics. For instance, he proves that two conics\index{conics} intersect in at most two points ``if their convexities are not in the same direction" (Proposition 30).\footnote{\emph{The Conics}, Book IV, \cite[p.\,172]{Apollonius-IV}.}  Proposition 35  of the same book concerns the tangency of conics with convexities in opposite directions.\footnote{\emph{The Conics}, Book IV, \cite[p.\,178]{Apollonius-IV}.}  
 Proposition 37 concerns the intersection of a hyperbola with another conic with convexities in opposite directions,\index{Apollonius of Perga (c. 262--c. 190 BCE)!\emph{Conics}}\footnote{\emph{The Conics}, Book IV, \cite[p.\,183]{Apollonius-IV}.} and there are several other examples where convexity is involved. We note that conics\index{conics} themselves are convex---they bound convex regions of the plane.

In the works of Archimedes,\index{Archimedes (c. 287--c. 212 BCE)} considerations on convexity\index{convexity} are not limited to conics but they concern arbitrary curves and surfaces. 
Right at the beginning of his treatise \emph{On the sphere and the cylinder}\index{Archimedes (c. 287--c. 212 BCE)!\emph{On the sphere and the cylinder}} \cite[p.\,2]{Heath-Archimedes}, Archimedes introduces the general notion of convexity. The first definition affirms the existence of ``bent lines in the plane which either lie wholly on the same side of the straight line joining their extremities, or have no part of them on the other side."\footnote{The quotations of Archimedes are from Heath's  translation \cite{Heath-Archimedes}.} Definition 2 is that of a concave curve:

\begin{quote}\small
I apply the term \emph{concave in the same direction} to a line such that, if any two points on it are taken, either all the straight lines connecting the points fall on the same side of the line, or some fall on one and the same side while others fall on the line itself, but none on the other side.
\end{quote}

Definitions 3 and 4 are the two-dimensional analogues of Definitions 1 and 2. Definition 3 concerns the existence of surfaces with boundary\footnote{``Terminated surfaces", in Heath's translation.} whose boundaries are contained in a plane and such that ``they will either be wholly on the same side of the plane containing their extremities,\footnote{In this and the next quotes, since we follow Heath's translation, we are using the word ``extremities", although the word ``boundary" would have been closer to what we intend in modern geometry.} or have no part of them on the other side."
Definition 4 is about concave surfaces, and it is an adaptation of the one concerning concave curves: 

\begin{quote}\small
I apply the term \emph{concave in the same direction} to surfaces such that, if any two points on them are taken, the straight lines connecting the points either all fall on the same side of the surface, or some fall on one and the same side of it while some fall upon it, but none on the other side.
\end{quote}

After the first definitions, Archimedes\index{Archimedes (c. 287--c. 212 BCE)} makes a few \emph{assumptions}, the first one being that among all lines having the same extremities, the straight line is the shortest. The second assumption is a comparison between the lengths of two concave curves in the plane having the same endpoints, with concavities in the same direction, and such that one is contained in the convex region bounded by the other curve and the line joining its endpoints:

\begin{quote}\small
Of other lines in a plane and having the same extremities, [any two] such are unequal whenever both are concave in the same direction and one of them is either wholly included between the other and the straight line which has the same extremities with it, or is partly included by, and is partly common with, the other; and that [line] which is included is the lesser.
\end{quote}

The third and fourth assumptions are analogues of the first two for what regards area instead of length. Assumption 3 says that among 
all surfaces that have the same extremities and such that these extremities are in a plane,  the plane is the least in area. Assumption 4 is more involved:

\begin{quote}\small
Of other surfaces with the same extremities, the extremities being in a plane, [any two] such are unequal whenever both are concave in the same direction and one surface is either wholly included between the other and the plane which has the same extremities with it, or is partly included by, and partly common with, the other; and that which is included is the lesser [in area].
\end{quote}

Almost all of Book I (44 propositions) of Archimedes'\index{Archimedes (c. 287--c. 212 BCE)}  treatise  \emph{On the sphere and the cylinder}\index{Archimedes (c. 287--c. 212 BCE)!\emph{On the sphere and the cylinder}} uses in some way or another the notion of convexity.\index{convexity} Several among these propositions are dedicated to inequalities concerning length and area under convexity assumptions. For instance, we find there  inequalities on the length of polygonal figures inscribed in convex figures.
Archimedes\index{Archimedes (c. 287--c. 212 BCE)} proves the crucial fact that the length of a convex curve is equal to the limit of polygonal paths approximating it, and similar propositions concerning area and volume, in particular for figures inscribed in or circumscribed to a circle or a sphere. 
His booklet \emph{On measurement of a circle}\index{Archimedes (c. 287--c. 212 BCE)!\emph{On measurement of a circle}} is on the same subject. His treatise \emph{On the equilibrium of planes}\index{Archimedes (c. 287--c. 212 BCE)!\emph{On the equilibrium of planes}}\index{Archimedes (c. 287--c. 212 BCE)!\emph{On the equilibrium of planes}} is another work in which convexity\index{convexity} is used in a fundamental way. Postulate 7 of Book I of that work says that ``in any figure whose perimeter is concave in one and the same direction, the center of gravity must be within the figure."

Another discovery of Archimedes\index{Archimedes (c. 287--c. 212 BCE)} involving the notion of convexity is his list of thirteen polyhedra that are now called semi-regular (Archimedean) convex polyhedra.  These are convex polyhedra whose faces are regular polygons of a non-necessarily unique type but admitting a symmetry group which is transitive on the set of vertices.\footnote{We are using modern terminology.} The faces of such a polyhedron may be a mixture of equilateral triangles and squares, or of equilateral triangles and regular pentagons, or of regular pentagons and regular hexagons, etc. It turns out that the Archimedean polyhedra are finite in number (there are essentially thirteen of them, if one excludes the regular ones).  Archimedes' work on this subject does not survive but Pappus,\index{Pappus of Alexandria (c. 290--c. 350)} in Book V of the \emph{Collection}\index{Pappus of Alexandria (c. 290--c. 350)!\emph{Collection}}, reports on this topic, and he says there  that Archimedes\index{Archimedes (c. 287--c. 212 BCE)} was the first to give the list of 13 semi-regular polyhedra \cite[p.\,272-273]{Pappus}.

Proclus (411--485)\index{Proclus of Lycia (412--485)} uses the notion of convexity\index{convexity} at several places of his \emph{Commentary on the first book of Euclid's\index{Euclid (c. 325--270 BCE)} Elements}\index{Proclus of Lycia (411--485)!\emph{In Euclidem}}; see e.g. his commentary on Definitions IV, VIII and XIX in which he discusses the angle made by two circles, depending on the relative convexities of the circles \cite[p.\,97, 115 and 141]{Proclus}.

Ptolemy,\index{Ptolemy (c. 100--c. 170)}  in Book I of his \emph{Almagest}\index{Ptolemy (c. 100--170)!\emph{Almagest}}, establishes a necessary and sufficient condition for a convex quadrilateral to be inscribed in a circle in terms of a single relation between the lengths of the sides and those of the diagonals of that quadrilateral.\footnote{Ptolemy's proof with the reference to Heiberg's edition is quoted in Heath's edition of Euclid \cite[Vol. 2 p. 225]{Euclid-Heath}.} The relation, known as \emph{Ptomely's relation}, has always been very useful in geometry.

\section{Mirrors and optics}
Convex and concave mirrors\index{mirror} are traditionally associated with imagination and phantasies, because they distort images. There are many examples of visual illusions and deceptions caused by convex and concave mirrors. Plato,\index{Plato (428--348 BCE)} in the \emph{Republic} (in particular, in the well-known cave passage),\footnote{Actually, in the cave passage (\cite{Plato-Republic}, Book VII, 514a-521d), not only images are distorted because the walls are not planar, but also one sees only shadows, apparent contours. Thom,\index{Thom, René (1923--2002)} in his \emph{Esquisse d'une sémiophysique} (\cite[p. 218]{Esquisse} of the English translation) sees there the mathematical problem of reconstructing figures from their apparent contours.}  uses this as an illustration of his view that reality is very different from sensible experience.  
According to him,\index{Plato (428--348 BCE)} reason, and especially mathematics, allows us to see the real and intelligible world of which otherwise we see only distorted shadows.

The Roman writer Pliny the elder\index{Pliny the elder (23--79 CE)} (1st c. CE) at several 
places of his \emph{Natural history} refers to the concept of concave or convex surface. In a passage on mirrors,\index{mirror} he writes \cite[Vol. VI, p. 126]{Pliny}:

\begin{quote}
Mirrors, too, have been invented to reflect monstrous 
forms; those, for instance, which have been consecrated in the 
Temple at Smyrna. This, however, all results from the configuration given to the metal; and it makes all the difference 
whether the surface has a concave form like the section of a 
drinking cup, or whether it is [convex] like a Thracian buckler; whether it is depressed in the middle or elevated; 
whether the surface has a direction transversely or obliquely; or whether it runs horizontally or vertically; the 
peculiar configuration of the surface which receives the shadows, causing them to undergo corresponding distortions: for, in 
fact, the image is nothing else but the shadow of the object 
collected upon the bright surface of the metal. 
\end{quote}

Regarding polyhedra receiving shadows, let me also mention the sundials\index{sundial} used in Greek antiquity that have the form of a convex surface (Figure \ref{Sundial}). The curve traced by the shadow of the sun has an interesting mathematical theory. The seventeenth-century mathematician Philippe de la Hire, in his treatise titled  \emph{Gnomonique ou l'art de tracer des cadrans ou horloges solaires sur toutes les surfaces, par différentes pratiques, avec les démonstrations géométriques de toutes les opérations} (Gnomonics, or the art of tracing sundials over all kind of surfaces by different methods, with geometrical proofs of all the operations) \cite{Hire}, conjectured that the theory of conic sections originated in the practical observations of sundials. Otto Neugebauer,\index{Neugebauer, Otto Eduard (1899--1990)} in his paper  \emph{The Astronomical Origin of the Theory of Conic Sections} \cite{Neugebauer}, made the same conjecture. This is also discussed in the article \cite{ACampo-Papadopoulos} in the present volume.

  \begin{center}
\includegraphics[width=0.8\linewidth]{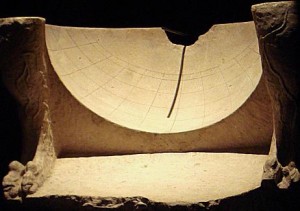} 
\vskip .1in \label{Sundial}
 {\small A Greek sundial with convex plate, from Ai Khanoum (Afghanistan), 3$^{\mathrm{rd}}$-2$^{\mathrm{nd}}$ c. BCE}
  \end{center}

With mirrors, we enter into the realm of optics,\index{optics} where convexity\index{convexity} is used in an essential way.  The propagation properties of light rays, including their reflection and their refraction properties on convex and concave mirrors,\index{mirror} were studied extensively by the Greek mathematicians.  Catoptrics,\index{catoptrics} the science of mirrors, was considered as a mathematical topic and was closely related to the theory of conic sections.\index{conics} (The Greek word katoptron, \tg{k'atoptron} means mirror.) Since ancient times, the study in this field involved, not only plane mirrors, but curved ones as well, concave or convex. 
A few books on catoptrics dating from Greek antiquity have reached us, some only in Arabic or Latin translation. There is a treatise with a possible attribution to Euclid,\index{Euclid (c. 325--270 BCE)} compiled  and amplified by Theon of Alexandria (IVth c. CE),\index{Theon of Alexandria (c. 335--405)} containing reflection laws for convex and concave mirrors; cf. \cite{Euclid-Ver-Eecke}.  In particular, the author studies there the position of the focus of a concave mirror, that is, the point where sun rays concentrate after reflection, so as to produce fire. 

The geometry of mirrors is related to conic sections. Book III of Apollonius'\index{Apollonius of Perga (c. 262--c. 190 BCE)}  \emph{Conics}\index{Apollonius of Perga (c. 262--c. 190 BCE)!\emph{Conics}} addresses the question of reflection properties of these curves. 
 A treatise by Heron of Alexandria\index{Heron of Alexandria (c. 10--70 CE)} which survives in the form of fragments is concerned with the laws of reflection on plane, concave and convex mirrors and their applications. In another treatise on catoptrics, attributed to Ptolemy\index{Ptolemy (c. 100--c. 170)} and whose 3${}^{\mathrm{rd}}$, 4${}^{\mathrm{th}}$ and 5${}^{\mathrm{th}}$  Books survive, the author studies the reflection properties on plane, spherical convex and spherical concave mirrors. A book titled \emph{Catoptrics}\index{catoptrics}\index{Archimedes (c. 287--c. 212 BCE)!\emph{Catoptrics}} by Archimedes\index{Archimedes (c. 287--c. 212 BCE)} does not survive but is quoted by later authors, notably by Theon of Alexandria in his \emph{Commentary on Ptolemy's Almagest} (I.3). One should also mention the work of Diocles (3${}^{\mathrm{rd}}$--2${}^{\mathrm{nd}}$ c. BCE) which survives in the form of fragments, citations by other authors, and translations and commentaries by Arabic mathematicians. Diocles\index{Diocles (240--180 BCE)} was a contemporary of Apollonius\index{Apollonius of Perga (c. 262--c. 190 BCE)}  and his work on optics\index{optics} is inseparable from the theory of conics.\index{conics}   To him is attributed the first investigation of the focal property of the parabola. (Heath in his edition of the \emph{Conics} \cite[p.\,114]{A-Heath-Conics} notes that Apollonius never used or mentioned the focus of a parabola.) Diocles studied mirrors having the form of pieces of spheres, paraboloids of revolution, ellipsoids of revolution, and other surfaces.  His work is edited by Rashed in \cite{Rashed-Catoptriciens}, a book containing critical editions of Arabic translations of  Greek texts on the theory of burning mirrors,\index{burning mirror}\index{mirror!burning} in particular those  by Diocles on elliptical and parabolic burning mirrors.\footnote{The Latin word \emph{focus} means fireplace, which led to the expression ``burning mirror."} The questions of finding the various shapes that a mirror can take in order to concentrate sun rays onto a point and produce fire at that point, and conversely, given a mirror, to find (possibly) a point where sun rays reflected on that mirror concentrate and produce fire, are recurrent in the Greek treatises on optics. There are proposition in Euclid's \emph{Catoptrics} \cite{Euclid-Ver-Eecke} dealing with burning mirrors. The introductory chapter of Lejeune's \emph{Recherches sur la catoptrique grecque d'après les sources antiques et médiévales} (Researches on Greek catoptrics following antique and medieval sources)  \cite{Lejeune} contains an interesting brief history of this subject.

 Optics is related to astronomy, in particular because of lenses.
 Convex lenses were already used in Greek antiquity to explore the heavens.  There is a famous passage of the  \emph{Life of Pythagoras}  written 
 by Iamblichus,\index{Iamblichus (245--325)}\index{Iamblichus (245--325)!\emph{Life of Pythagoras}} the Syrian neo-Pythagorean mathematician of the third century CE, in which the author recounts that Pythagoras, at the moment he made his famous discovery of the relation between ratios  of integers and musical intervals, was pondering on the necessity of finding a device which would be useful for the ear in the same manner as the dioptre is useful for the sight \cite[p.\,62]{Iamblichus}. Dioptres are a kind of glasses used to observe the celestial bodies.

\section{Billiards in convex domains}
A question raised by Ptolemy\index{Ptolemy (c. 100--c. 170)} is known since the Renaissance as \emph{Alhazen's problem}.\footnote{\label{f:Hatham}The name refers to  Ibn al-Haytham,\index{Ibn al-Haytham (965--1040)} the Arab scholar from the Middle Ages known in the Latin world as Alhazen, a deformation of the name ``Al-Haytham." Ibn  al-Haytham is especially famous for his treatise on Optics (\emph{Kit\=ab al-man\a=azir}), in seven books (about 1400 pages long), which was translated into Latin at the beginning of the thirteenth century, and which was influential on Johannes Kepler, Galileo Galilei, Christiaan Huygens and Ren\'e Descartes, among others. 
An important part of what survives from his work in geometry and optics\index{optics} was translated and edited by Rashed \cite{Rashed-IV-1, Rashed-V}. Ibn al-Haytham\index{Ibn al-Haytham (965--1040)}   is the author of an ``intromission"\index{vision!intromission theory} theory of vision saying that it is the result of light rays penetrating our eyes, contradicting the theories held by Euclid\index{Euclid (c. 325--270 BCE)} and Ptolemy who considered, on the contrary, that vision is the result of light rays emanating from the eye (``extramission" theory).\index{vision!extramission theory} It is possible though that Euclid, as a mathematician, adhered to the theory where visual perception is caused by light rays traveling along straight lines emitted from the eye that strike the objects seen, in oder to develop his mathematical theory of optics as an application of Euclidean geometry. This also explains the fact that Euclid's optics does not include any physiological theory of vision, nor any physical theory of colors, etc. Needless to say, besides this rough classification into an intromission theory and an extramission theory of light, there is a large amount of highly sophisticated and complex theories of vision and of light that were developed by Greek authors, which were related to the various philosophical schools of thought, and at the same time to the mathematical theories that were being developed.}
This problem, in its generalized form, concerns reflection in a convex mirror, and, in the modern terminology, it can be regarded as a problem concerning trajectories in a convex billiard table. In its original form, the problem asks for the following: given a circle and two points that both lie outside or inside this circle, to construct a point on the circle such that two lines joining the given points to that point on the circle make equal angles with the normal to the circle at the constructed point. Ibn al-Haytham,\footnote{See Footnote \ref{f:Hatham}.} in his \emph{Kit\=ab al-Man\=azir}, found a geometrical solution of the problem using conic sections and asked for an algebraic solution. He referred to Ptolemy while writing on his problem, and in fact, a large part of his work on optics was motivated by Ptolemy's work on this subject, which he criticized at several points.  A. I. Sabra published an edition of Ibn al-Haytham's\index{Ibn al-Haytham (965--1040)} \emph{Optics} \cite{Sabra-Kitab}, and wrote a paper containing an account of six lemmas used by Ibn al-Haytham in his work on the problem on the problem \cite{Sabra, Simon, Smith, Smith2001}.

A number of prominent mathematicians worked on Alhazen's problem. We mention in particular Christiaan Huygens,\index{Huygens, Christiaan (1629--1695)} who wrote several articles and notes on it; they are published in Volumes XX and XXII of his \emph{Complete Works} edition \cite{OC}.  J. A. Vollgraff, the editor of Volume XXII of these works, writes on p. 647:   ``At the beginning of 1669, we can see Huygens absorbed by mathematics. He was busy with Alhazen's problem. This is one of the problems of which he always strived to find, using conic sections, the most elegant solution."\footnote{Au commencement de 1669 nous voyons Huygens absorbé par la mathématique. Il s'occupa du problème d'Alhazen. C'est là un des problèmes dont il a toujours eu l'ambition de trouver, par les sections coniques, la solution la plus élégante.} Volume XX of  the \emph{Complete Works} contains a text read by Huygens in 1669 or 1670 at the Royal Academy of Sciences of Paris on this problem titled \emph{Problema Alhaseni} (p. 265). There is also a note titled  \emph{Construction d'un Problème d'Optique, qui est la XXXIX$^{\mathrm{e}}$ Proposition du Livre
V d' Alhazen, et la XXII$^{\mathrm{e}}$ du Livre VI de Vitellion} (Construction of a problem on optics, which is Proposition XXXIX of Book V of Alhazen and Proposition XXII of Book VI of Vitellion)\footnote{Vitellion is the name of a thirteenth-century mathematician who edited works of Alhazen on optics.} and a note (p. 330) titled \emph{Problema Alhazeni ad inveniendum in superficie speculi sphaerici punctum reflexionis} (Alhazen's problem on finding the reflection point on the surface of a spherical mirror), in the same volume. Vol. XXII of the \emph{Complete Works} \cite{OC} contains an article dating from 1673  titled \emph{Constructio et demonstratio ad omnes casus Problematis Alhazeni de
puncto reflexionis} (Construction and proof of all the cases of  Alhazen's reflection point problem). There are other pieces related to Alhazen's problem in Huygens' complete works. Among the other mathematicians who worked on this problem, we mention the Marquis de l'H\^opital and Isaac Barrow. Leonardo da Vinci conceived a mechanical device to solve the problem. Talking about Leonardo, let us note that he also conceived devices to draw conics; cf. \cite{Pedretti}.

    \section{Aristotle}
 
 Aristotle\index{Aristotle (384--322 BCE)} mentions convex and concave surfaces at several places in his writings, usually in his explanations of metaphysical ideas, and generally as an illustration of the fact that the same object may have two very different appearances, depending on the way one looks at it; like the circle, which may seem concave or convex, depending on the side from which one sees it.

In the \emph{Nicomachean Ethics}\index{Aristotle (384--322 BCE)!\emph{Nicomachean ethics}} \cite{A-Nicomachean},\footnote{\emph{Nicomachean Ethics} \cite{A-Nicomachean}, 1102a-30.} talking about the soul, Aristotle\index{Aristotle (384--322 BCE)} discusses the fact that it has a part which is rational and a part which is irrational.  He asks whether this distinction into two parts is comparable to the distinction between the parts of a body or of anything divisible into parts, or whether these two parts are 
``by nature inseparable, like the convex and concave parts in the circumference of a circle." The latter response is, according to him, the correct one, because the soul is one, and the fact that it has rational and irrational behaviors are different phases of the same thing.

Convexity\index{convexity} is also mentioned in the \emph{Meteorology}\index{Aristotle (384--322 BCE)!\emph{Meteorology}} \cite{A-Meteorology}. Here, Aristotle\index{Aristotle (384--322 BCE)} explains the origin of rivers and springs. He writes:\footnote{\emph{Meteorology}  \cite{A-Meteorology}, 350a10.} 
\begin{quote}\small
For mountains and high ground, suspended over the country like a saturated sponge, make the water ooze out and trickle together in minute quantities but in many places. They receive a great deal of water falling as rain (for it makes no difference whether a spongy receptacle is concave and turned up or convex and turned down: in either case it will contain the same volume of matter) and, they also cool the vapour that rises and condense it back into water.
\end{quote}

 Chapter 9 of Aristotle's\index{Aristotle (384--322 BCE)} \emph{Physics}\index{Aristotle (384--322 BCE)!\emph{Physics}} \cite{A-Physics} is concerned with the existence of void, and it is an occasion for the Philosopher to discuss actuality and potentiality, in various instances. Convexity\index{convexity} and concavity are again used as metaphorical entities. He writes:\footnote{\emph{Physics} \cite{A-Physics}, 217a30-b5.} 
 \begin{quote}\small
 For as the same matter becomes hot from being cold, and cold from
being hot, because it was potentially both, so too from hot it can
become more hot, though nothing in the matter has become hot that
was not hot when the thing was less hot; just as, if the arc or curvature
of a greater circle becomes that of a smaller---whether it remains
the same or becomes a different curve---convexity has not come to exist
in anything that was not convex but straight.
\end{quote}
 Thus, he says, something which becomes cold and hot was potentially cold or hot, like a thing which is convex: it may become more convex, because it was potentially convex, but it cannot become straight. Chapter 13 of the same treatise is concerned with the meaning of different words related to \emph{time}. Aristotle\index{Aristotle (384--322 BCE)} writes\index{Aristotle (384--322 BCE)!\emph{Physics}}:\footnote{\emph{Physics} \cite{A-Physics},  222b1.}  
 \begin{quote}\small
 Since the ``now" is an end and a beginning of time, not of the same
time however, but the end of that which is past and the beginning
of that which is to come, it follows that, as the circle has its convexity\index{convexity}
and its concavity, in a sense, in the same thing, so time is always
at a beginning and at an end.
\end{quote}
 Thus, again like a circle, depending on the side from which one looks at it, may seem concave or convex, the term ``now", depending on the side from which we look at it, may be the beginning or the end of time.

 Besides talking about the convex and the concave as two sides of the same thing, Aristotle\index{Aristotle (384--322 BCE)} liked to give the example of the convex and the concave, while he talked about opposites. In the \emph{Mechanical problems}\index{Aristotle (384--322 BCE)!\emph{Mechanical problems}} \cite{A-Mechanica}, talking about the lever, and the fact that with this device, a small force can move a large weight, he writes:\footnote{\emph{Mechanical problems} \cite{A-Mechanica}, 847b.} 
 \begin{quote}\small
 The original cause of all such phenomena is the circle; and this is natural, for it is in no way strange that something remarkable should result from something more remarkable, and the most remarkable fact is the combination of opposites with each other. The circle is made up of such opposites, for to begin with it is composed both of the moving and of the stationary, which are by nature opposite to each other.[\ldots] an opposition of the kind appears, the concave and the convex. These differ from each other in the same way as the great and small; for the mean between these latter is the equal, and between the former is the straight line.
 \end{quote}
  Aristotle's\index{Aristotle (384--322 BCE)} fascination for the circle is always present in his works. In another passage of the same treatise, he writes\index{Aristotle (384--322 BCE)!\emph{Mechanical problems}}:\footnote{\emph{Mechanical problems} \cite{A-Mechanica}, 848a.} 
 \begin{quote}\small
 [the concave and the convex] before they could pass to either of the extremes, so also the line must become straight either when it changes from convex to concave, or by the reverse process becomes a convex curve. This, then, is one peculiarity of the circle, and a second is that it moves simultaneously in opposite directions; for it moves simultaneously forwards and backwards, and the radius which describes it behaves in the same way; for from whatever point it begins, it returns again to the same point; and as it moves continuously the last point again becomes the first in such a way that it is evidently changed from its first position.
 \end{quote}

 Book V of Aristotle's\index{Aristotle (384--322 BCE)} \emph{Problems}\index{Aristotle (384--322 BCE)!\emph{Problems}} \cite{A-Problems}  is titled \emph{Problems connected with fatigue}. Problem 11 of that book asks: ``Why is it more fatiguing to lie down on a flat than on 
a concave surface? Is it for the same reason that it is 
more fatiguing to lie on a convex than on a flat surface?" As usual in Aristotle's\index{Aristotle (384--322 BCE)} \emph{Problems}, the question is followed by comments and partial answers, some of them may be due to Aristotle,\index{Aristotle (384--322 BCE)} and others presumably written by students of the Peripatetic school. The comments on this problem include a discussion on the pressure exerted on a convex line, saying that it is greater  than that exerted on a straight or concave line. They start with: 
\begin{quote}\small For the weight being concentrated in one place in the sitting 
or reclining position causes pain owing to the pressure. This 
is more the case on a convex than on a straight surface, and 
more on a straight than on a concave; for our body assumes 
curved rather than straight lines, and in such circumstances 
concave surfaces give more points of contact than flat 
 surfaces. For this reason also couches and seats which 
yield to pressure are less fatiguing than those which do not 
do so.
\end{quote}
 Book XXXI of the same work is titled  \emph{Problems connected with the eyes}, and Problem 25  of that book involves a discussion of convexity\index{convexity} in relation with vision. The question is:  ``Why is it that though both a short-sighted and an old man 
are affected by weakness of the eyes, the former places an 
object, if he wishes to see it, near the eye, while the latter 
holds it at a distance?" In the comments, we read: 
\begin{quote}\small The short-sighted man can see the object but 
cannot proceed to distinguish which parts of the thing
at which he is looking are concave and which convex, 
 but he is deceived on these points. Now concavity and 
convexity\index{convexity} are distinguished by means of the light which 
they reflect; so at a distance the short-sighted man cannot 
discern how the light falls on the object seen; but near at 
hand the incidence of light can be more easily perceived.
\end{quote}

The treatise \emph{On the Gait of Animals} is a major biological treatise of Aristotle, and it concerns motion and the comparison of the various ways of motion for animals  (including human beings). 
In Chapter 1, Aristotle\index{Aristotle (384--322 BCE)} writes\index{Aristotle (384--322 BCE)!\emph{On the gait of animals}}:\footnote{\emph{On the Gait of Animals} \cite{A-Gait} 704a15.} 
\begin{quote}\small Why do man and bird, though both bipeds, have an opposite curvature of the legs? For man bends his legs convexly, a bird has his bent concavely; again, man bends his arms and legs in opposite directions, for he has his arms bent convexly, but his legs concavely. And a viviparous quadruped bends his limbs in opposite directions to a man's, and in opposite directions to one another; for he has his forelegs bent convexly, his hind legs concavely.
\end{quote}
 In Chapter 13 of the same treatise, we read:\footnote{\emph{On the Gait of Animals} \cite{A-Gait}, 712a.}  
 \begin{quote}\small There are four modes of flexion if we take the combinations in pairs. Fore and hind may bend either both backwards, as the figures marked A, or in the opposite way both forwards, as in B, or in converse ways and not in the same direction, as in C where the fore bend forwards and the hind bend backwards, or as in D, the opposite way to C, where the convexities are turned towards one another and the concavities outwards.
 \end{quote}

\section{Architecture} The Parthenon, Erechteum and Theseum columns, and more generally, Doric columns, are not straight but convex. Several explanations for this fact have been given, but none of them is definitive. F. C. Penrose published a book titled  \emph{An investigation of the principles of Athenian architecture; or, the results of a
survey conducted chiefly with reference to the optical refinements exhibited in the construction
of the ancient buildings at Athens} \cite{Penrose}, whose object he describes as (p. 22) 
\begin{quote}\small
the investigation of various delicate curves, which form the principal architectural lines of certain of the Greek buildings of the best period; which lines, in ordinary architecture, are (or are intended to be) straight. In the course of our inquiries we shall perhaps be enabled in some degree to extend and correct our views of the geometry and mathematics of the ancients, by establishing the nature of the curves employed [\ldots] The most important curves in point and extent are those which form the horizontal lines of the buildings where they occur; such as the edges of the steps and the lines of the entablature, which are usually understood to be straight level lines, but in the steps of the Parthenon and some other of the best examples of Greek doric, are convex curves [\ldots].
\end{quote}

 For instance, the columns of the Parthenon are ``in wonderful agreement at all points" with a piece of a parabola \cite[p.\,41]{Penrose}.

\begin{center}
\includegraphics[width=\linewidth]{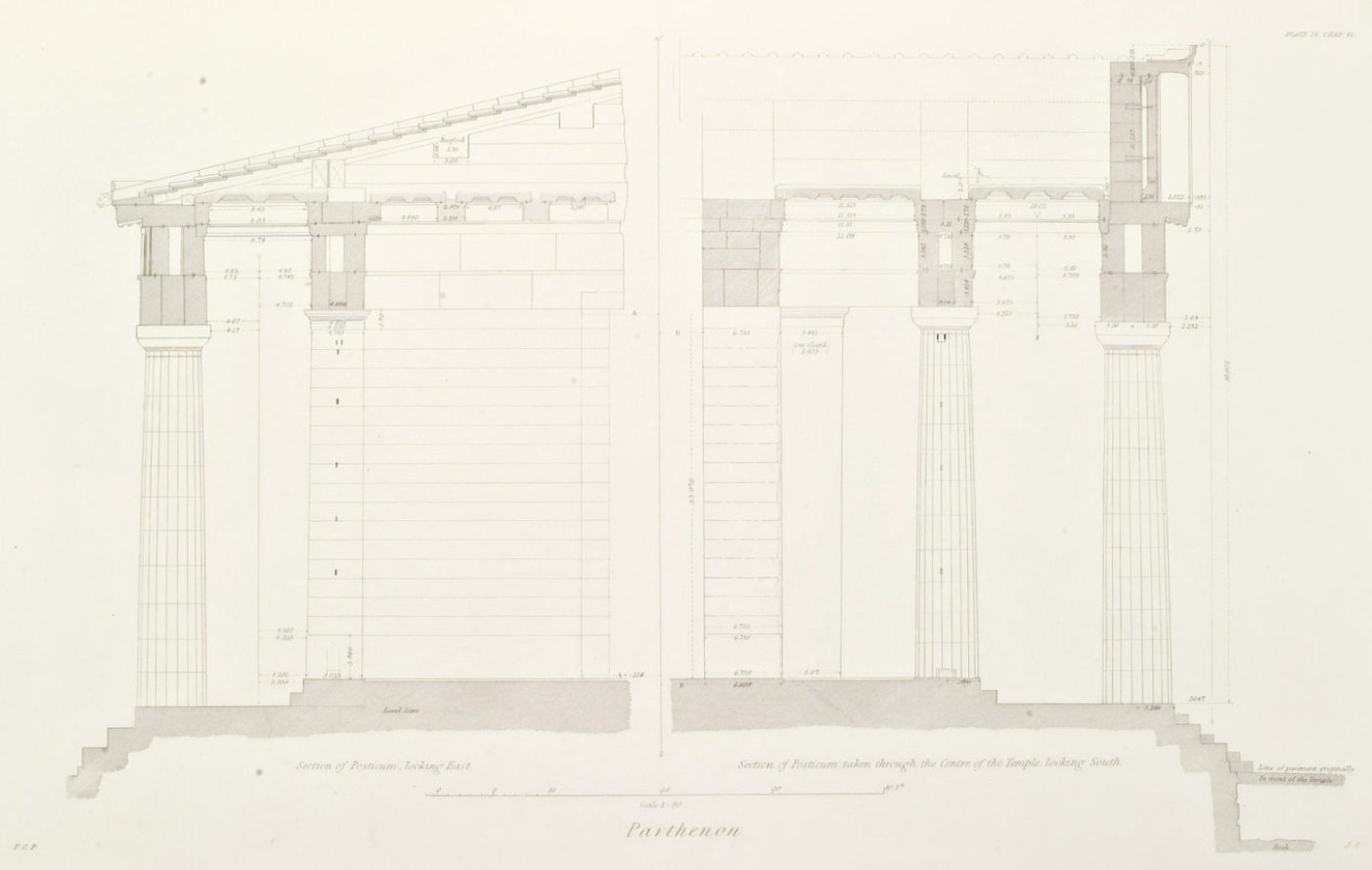} 
\vskip .1in \label{Parthenon}
 {\small A drawing of the Parthenon, from F. C. Penrose's \emph{Investigation of the Principles of Athenian Architecture}.} 
  \end{center}
  \vskip .3in 

Penrose says that the first mention of the curvature properties in Greek architecture was made by the Roman historian of architecture Vitruvius.\index{Vitruvius (Marcus Vitruvius Pollio) (I$^{\mathrm{st}}$ c. BCE)} Referring to him again on p. 39 of his essay \cite{Penrose}, Penrose writes that this phenomenon, called \emph{entasis} (from a Greek word meaning to stretch a line, or to bend a bow), is the 
\begin{quote}\small
well-known increment or swelling given to a column in the middle parts of the shaft for the purpose of correcting a disagreeable optical illusion, which is found to give an attenuated appearance to columns formed with straight sides, and to cause their outlines to seem concave instead of straight. The fact is almost universally recognized by attentive observers, though it may be difficult to assign a conclusive reason why it should be so.
\end{quote} Another  possible explanation which Penrose gives is ``simply an imitation of the practice of Nature in giving almost invariably a convex outline of the limbs of animals" (p. 116) and ``of the trunks and branches of trees" (p. 105).

This leads us to the question of architecture and art as an imitation of Nature, whose lines are seldom straight, and sometimes intricately curved. Curved are also the roads of poetic creation. The first lines of Canto I of Dante's\index{Dante Alighieri (1265--1321)} \emph{Inferno} read: 
\begin{quote} 
Midway upon the journey of our life
\\
I found myself within a forest dark,
\\
For the straightforward pathway had been lost.\footnote{ \emph{Nel mezzo del cammin di nostra vita 
 \\ \hskip .2in 
mi ritrovai per una selva oscura 
\\ \hskip .2in 
ché la diritta via era smarrita.}}
\end{quote}

The roads of mathematical discovery are  even more curved; they are twisted, and very lengthy.

%
%
%
%

  \printindex

\end{document}